\documentstyle{amsppt}
\magnification\magstep1
\NoRunningHeads
\pageheight{9 truein}
\pagewidth{6.3 truein}
\baselineskip=15pt
\catcode`@=11
\def\logo@{\relax}
\catcode`@=\active
\def\af #1.{\Bbb A^{#1}}
\def\au#1.{\operatorname {Aut}\,(#1)}

\def\ring#1.{\Cal O_{#1}}
\def\alb#1.{\ring .(#1)}
\def\pr #1.{\Bbb P^{#1}}
\def\prr #1.{{\Bbb P}(#1)}

\def\pic#1.{\operatorname {Pic}\,(#1)}

\def\mg{\overline {\Cal M}_g}

\def\picG#1.{{\operatorname {Pic}}^G(#1)}
\def\psl#1.{\operatorname {PSL}_{#1}}
\def\Sl#1.{\operatorname {SL}_{#1}}
\def\mov#1.{\operatorname {Mov}(#1)}
\def\nef#1.{\operatorname {Nef}(#1)}
\def\nd#1.{N^1(#1)}
\def\nsd#1.{\nd #1.}
\def\nsdg#1.{N^1_G(#1)}
\def\nc#1.{N_1(#1)}
\def\mc#1.{\overline{NE}_1(#1)}
\def\ec#1.{\overline{NE}^1(#1)}
\def\mcd#1.{\ec #1.}
\def\deq{:=}
\def\ex#1.{\operatorname{ex}(#1)}
\def\sym#1.#2.{\operatorname{sym}_{#1}(#2)}
\def\spec#1.{\operatorname{spec}(#1)}
\def\Hom#1.#2.{\operatorname{Hom}(#1,#2)}
\def\bbbq{\Bbb Q}
\def\bbba{\Bbb A}
\def\gm{{\Bbb G}_m}
\def\cox#1.{\operatorname{Cox}(#1)}

\def\bbbn{\Bbb N}
\def\bbbz{\Bbb Z}
\def\bbbc{\Bbb C}
\def\bbbp{\Bbb P}

\def\sr#1.#2.{\operatorname{R}(#1,#2)}
\def\proj#1.{\operatorname{Proj}(#1)}
\def\ac #1.{\operatorname{Ample}(#1)}
\def\gac #1.{C^G(#1.)}

\def\tm #1.{\operatorname{taut}(#1)}
\def\picgq#1.{\operatorname {Pic}^G(#1)_{\bbbq}}
\def\picq#1.{\operatorname {Pic}(#1)_{\bbbq}}
\def\picg#1.{\picG {#1}.}
\def\mon{\overline{M}_{0,n}}
\def\tmon{\widetilde{M}_{0,n}}
\def\fm #1.#2.{#1[#2]}
\def\fpone{{\pr 1.}[n]}
\def\ponen{(\pr 1.)^{\times n}}
\NoBlackBoxes

\topmatter
\title Mori Dream Spaces and GIT \endtitle
\author Yi Hu and Se\'an Keel \endauthor
\endtopmatter

$\;\;\;\;\;\;\;\;\;\;\;\;\;\;\;\;\;\;\;\;\;\;\;\;\;\;\;\;\;\;\;\;\;\;\;\;\;\;\;\;\;\;\;\;\;\;\;\;\;\;\;\;\;\;\;\;\;\;\;\;\;\;\;\;\;\;\;\;\;\;\;\;\;\;\;\;\;\;
\;\;\;\;\;\;\;\;\;\;\;\;\;\;${\it Dedicated to Bill Fulton }

An important advance 
 in algebraic geometry  in the last ten years 
is the theory of variation of geometric invariant theory quotient
(VGIT), see \cite{BP}, \cite{H}, \cite{DH},
\cite{T}. 
Several authors, have observed that
VGIT has implications for birational geometry, e.g. it gives
natural examples of Mori flips and contractions,  \cite{R2}, \cite{DH}, \cite{T}.
In this paper we observe that the connection is quite fundamental -- Mori theory is,
at an almost
tautological level, an instance of VGIT, see (2.14). Here are more details:

Given a projective variety $X$ a natural problem
is to understand the collection of all morphisms (with 
connected fibres) from $X$
to other projective varieties. Ideally one would like to decompose
each map into simple steps, and parameterize the possibilities,
both for the maps, and for the factorizations of each map. 
An important insight, principally
of Reid and Mori, is that the picture is often simplified if one allows
in addition to morphisms, {\it small modifications}, i.e. rational
maps that are isomorphisms in codimension one. With this extension
a natural framework is the category of rational contractions. In
many cases there is a nice combinatorial parameterization, given
by a decomposition of a convex polyhedral cone, the cone of effective
divisors $\ec X.$, into convex polyhedral chambers, which we
call Mori chambers. 
Instances of this structure have been
studied in various circumstances:
The existence of such a parameterizing decomposition
for Calabi-Yau manifolds was conjectured by Morrison, motivated
by ideas in mirror symmetry \cite{M}. The conjecture was
proven in dimension three by Kawamata, \cite{Ka}. Oda and Park
study the decomposition for toric varieties, motivated by
questions in combinatorics, \cite{OP}.   Shokurov studies such
a decomposition for parameterizing log minimal models, \cite{Sh}.
In Geometric Invariant Theory there is a similar combinatorial structure,
a decomposition of the $G$-ample cone into GIT chambers parameterizing
GIT quotients, \cite{DH}. The main observation of this paper is that whenever
a good Mori chamber decomposition exists, it is in a natural way
a GIT decomposition. 

The main goal of this paper is to study varieties $X$ with a good
Mori chamber decomposition (see \S 1 for the meaning of {\it good}).
We call such varieties
{\it Mori Dream Spaces}.
There turn out to be many examples, including
quasi-smooth projective toric (or more generally, spherical)
varieties, many GIT quotients, and log Fano $3$-folds.
We will
show that a Mori dream space is in a natural way a 
GIT quotients of affine variety 
by a torus 
in a manner generalizing Cox's construction of toric varieties as
quotients of affine space, \cite{C}. Via the quotient description,
the chamber decomposition of the cone of divisors is naturally 
identified with the decomposition of the $G$-ample cone from VGIT,
See (2.9). In particular {\bf every} rational contraction of a 
Mori dream space comes from GIT, and all possible factorizations
of a rational contraction (into other contractions) can be read off
from the chamber decomposition. See (2.3), (2.9) and (2.11).

Content overview: In \S 1 we define Mori chambers and Mori
dream spaces. The main theorems are proven in \S 2. In 
\S 3 we note connections with a question of Fulton about $\mon.$

Thanks: We would like to thank A. Vistoli, M. Thaddeus, J. Tate,
C. Teleman, and F. Rodriguez-Villegas for helpful discussions.
While working on this paper S.K. received partial support from
the NSF and NSA

\subhead \S 1 Mori equivalence \endsubhead
Throughout the paper $\nd X.$ indicates the Neron-Severi group
of divisors, with rational coefficients.

We begin with a few definitions:

\definition{1.0 Definitions-Lemma} Let $f: X \dasharrow Y$
be a rational map between normal projective varieties. 
Let 
$(p,q): W \rightarrow X \times Y$ be a resolution of $f$ with, $W$ projective,
and $p$ birational. We say $f$ has connected fibres if $q$ does. If
$f$ is birational, we call it a {\it birational contraction} if every
$p$ exceptional divisor is $q$ exceptional. 
For a $\bbbq$-Cartier divisor $D \subset Y$, 
$f^*(D)$ is defined to be $p_*(q^*(D))$. All of these are independent
of the resolution.
\enddefinition

Warning: $f^*$ for rational maps is not in general
functorial.

It is useful to generalize the notion of rational contraction to the
non-birational case. Intuitively this should be a composition of a 
small modification (see (1.8) below) and a morphism. Our definition is
different --we do not want to assume at the outset the existence of
small modifications  but in the cases we consider it will be equivalent,
see (1.11) below.

\definition{1.1 Definition} Notation as in (1.0). An effective divisor
$E$ on $W$ is called $q$-fixed if no effective Cartier divisor whose support
is contained in the support of $E$ is $q$-moving (see \cite{Ka}),
i.e. for every such divisor $D$ the natural map
$$
\ring Y. \rightarrow q_*(\ring .(D))
$$
is an isomorphism.
$f$ is called a {\it contraction} if every $p$ exceptional divisor
is $q$-fixed. An effective divisor $E \subset X$
is called $f$-fixed if any effective divisor of $W$ supported on the 
union of the strict transform of $E$ with the exceptional divisor of $p$
is $q$-fixed.
\enddefinition

One checks easily that for birational maps a divisor is fixed iff
it is exceptional.

\definition{1.2 Definition} For a line bundle $L$ on a scheme $X$
the section ring is the graded ring
$$
\sr X. L. \deq \bigoplus_{n \in \bbbn} H^0(X,L^{\otimes n}).
$$
\enddefinition

We will often mix the notation of divisors and line bundles, e.g.
writing $H^0(X,D)$ for $H^0(X,\ring.(D))$ for a divisor $D$.
We recall that the moving cone $\mov X. \subset \ec X.$ is 
the collection of (numerical
classes of) divisors with no stable base components.

If $\sr X. D.$ is finitely generated, and $D$ is effective, 
then there is an induced rational
map 
$$
f_D : X \dasharrow \proj {\sr X. D.}.
$$
which is regular outside the stable base locus of $|nD|$.

\definition{1.3 Definition-Lemma: Mori Equivalence} Let $D_1$ and $D_2$ be
two $\bbbq$-Cartier divisors on $X$ with finitely generated
section rings.
Then we say $D_1$
and $D_2$ are {\it Mori equivalent} if the rational maps $f_{D_i}$ have
the same Stein factorization 
i.e. there is an isomorphism between their images which makes
the obvious triangular diagram commutative. This occurs iff the rational
maps $f_{mD_i}$ are the same for some $m > 0$. 
\enddefinition

\definition{1.4. Definition} Let $X$ be a projective variety such that
$\sr X. L.$ is finitely generated for all line bundles $L$ and
$\picq X. = \nd X.$.  By a 
{\it Mori chamber} of $\nd X.$ we mean the closure of an equivalence class
whose interior is open in $\nd X.$.
\enddefinition

Contractions and finite generation turn out to be closely
related.

\proclaim{1.5 Lemma} Let $R = \bigoplus_{n \in \bbbn} R_n$ be  an
$\bbbn$-graded ring, finitely generated as an algebra over $R_0$.
Then for some $m>0$ the natural map

$$
\sym k. R_m. \rightarrow R_{km}
$$
is surjective for all $k > 0$. 
\endproclaim
\demo{Proof} Let $Y \deq \proj R.$. Then for some $m > 0$,
$H^0(Y,\ring Y.(km)) = R_{km}$ for all $k>0$. The result follows. \qed \enddemo

\proclaim{1.6 Lemma} If a divisor $D$ has a finitely generated section ring, then,
after replacing $D$ by a positive multiple, 
$f_D$ is a contracting rational map and
$D = f_D^*(\ring.(1)) + E$, for some $f_D$ fixed effective divisor
$E$. Conversely, if $f: X \dasharrow Y$
is a contracting rational map and $D = f_D^*(A) + E$ for $A$ ample
on $Y$, and $E$ fixed by $f$, then 
$D$ has a finitely generated section ring, and 
$f = f_{mD}$ for some $m > 0$. 
\endproclaim
\demo{Proof}
Suppose $\sr X. D.$ is finitely generated. Let $D = M + F$ be the canonical
decomposition of $D$ into its moving and fixed components. By (1.5) after
replacing 
$D$ by a multiple 
$$
\sym k. {H^0(X,M)}. \rightarrow H^0(X,kM)
$$
 is surjective 
and 
$$
H^0(X,kM) \rightarrow H^0(X,kM + rF) \tag{1.6.1}
$$
is an isomorphism for any $k,r > 0$.  By
passing to the blowup of the scheme-theoretic base locus of
$|M|$ we may assume $f$ is regular, and $M = f^*(A)$ for
some ample $A$ on $Y$. Now $F$ is $f$-fixed by (1.6.1). 

Now consider the converse, with notation as in the statement.
Let $p,q,W$ be a
resolution as in (1.0). By negativity of contraction,
\cite{Koll\'aretal,2.19} $p^*f^*(A) = q^*(A) + E'$ where $E'$ is
$p$ exceptional and effective. Thus $p^*(D) = q^*(A) + E''$
where $E''$ is effective and $q$-fixed. Thus
$\sr X. D. = \sr W. p^*(D). = \sr Y. A.$ is finitely generated. 
We can check $f = f_{mD}$ after throwing away the base locus of $D$,
where the equality is familiar. \qed \enddemo 

\proclaim{1.7 Lemma} Let $f: X \dasharrow Y$ and $g: X \dasharrow Z$
be birational contractions. Suppose $f^*(A) + E = g^*(B) + F$ for
$A$ ample, $B$ nef, $E$ $f$-exceptional and $F$ $g$ -exceptional.
then $ f \circ g^{-1}: Z \rightarrow Y$ is regular.
\endproclaim
\demo{Proof} By negativity of contractions we can pass to a resolution
and assume $f$ and $g$ are regular. Now by negativity of contraction
$E = F$, and so $f^*(A) = g^*(B)$. The result follows from the
rigidity lemma, see e.g. \cite{Ke,1.0}. \qed \enddemo

\definition{1.8 Definition} By a {\it small $\bbbq$-factorial
modification (SQM)} of a projective variety $X$ we mean a
contracting birational map $f : X \dasharrow X'$ with $X'$
projective and $\bbbq$-factorial, such that $f$ is an 
isomorphism in codimension one. \enddefinition

The most important examples of SQMs are flips.

\definition{1.9 Definition} Let 
$q: X \rightarrow Y$ a small birational morphism, and let $D$ be
a $\bbbq$-Cartier divisor such that $-D$ is $q$-ample. By a $D$-flip
of $q$ we mean a small birational morphism
$q':X' \rightarrow Y$ such that 
the strict transform
of $D$ on $X'$ is $\bbbq$-Cartier and $q'$ ample. We say the flip is
of relative Picard number one if $q$ and $q'$ are of relative Picard number one.
\enddefinition

The $D$-flip if it exists is unique. In the relative Picard number one
case it is independent of $D$. See e.g. \cite{KM}.

\definition{1.10 Definition: Mori Dream Space} We will
call a normal projective variety $X$ a Mori Dream Space 
provided the following hold:
\roster
\item $X$ is $\bbbq$-factorial and $\picq X. = \nd X.$
\item $\nef X.$ is the affine hull of finitely many semi-ample
line bundles. 
\item There is a finite collection of SQMs $f_i: X \dasharrow X_i$
such that each $X_i$ satisfies (2) and $\mov X.$ is the union
of the $f_i^*(\nef X_i.)$.
\endroster
\enddefinition

\proclaim{1.11 Proposition} Let $X$ be a Mori dream space. The following
hold:

\roster
\item Mori's
program can be carried out for any divisor on $X$, i.e. 
the necessary contractions and flips exist, any sequence terminates, and if
at some point the divisor becomes nef then at that point it becomes
semi-ample. 
\item The $f_i$ of (1.10) are the only SQMs of $X$. $X_i$, $X_j$ in adjacent
chambers are related by a flip. 
$\ec X.$ is the affine hull of finitely
many effective divisors. 
There are finitely many birational contractions
$g_i: X \dasharrow Y_i$ with $Y_i$ Mori dream spaces, such that
$$
\ec X. = \bigcup_{i} g_i^*(\nef Y_i.) \times \ex g_i.
$$
is a decomposition of $\ec X.$ into closed convex chambers, with
disjoint interiors. The cones $g_i^*(\nef Y_i.) \times \ex g_i.$ are precisely
the Mori chambers of $\ec X.$. They are in one to one correspondence
with birational contractions of $X$ with $\bbbq$-factorial image. 
\item The chambers $f_i^*(\nef X_i.)$ together with their faces
gives a fan, with support $\mov X.$. The cones in
the fan are in one to one correspondence with
contracting rational maps
$g: X \dasharrow Y$ with $Y$ normal and projective via
$$
[g: X \dasharrow Y ] \rightarrow [g^*(\nef Y.) \subset \mov X.] .
$$
\endroster
Let $D$ be an effective  divisor on $X$
\roster
\item[4] $\sr X. D.$ is finitely generated. 
\item After replacing $D$ by a multiple,
the canonical decomposition $D = M + F$ into moving and fixed part
has the following properties: There is a Mori
chamber containg $D$, so that if 
$g_i: X \dasharrow Y_i$ is the corresponding birational
contraction of (2), then $F$ has support the exceptional
locus of $g_i$ and $M$ is the pullback of a semi-ample
line bundle on $Y_i$. 
\endroster
\endproclaim
\demo{Proof} These all follow from the definition and purely formal
properties of Mori's program. Here is a sketch of the proof.

Note if $f: X \rightarrow Y$
is a small birational morphism, then $f^*(A)$ for $A$ ample 
is in the interior of $\mov X.$. Thus from (1.10.3)
all the small contractions of any $X_i$ have a flip, given by
another $X_j$. Now let $D$ be a divisor. If it is nef it is
semi-ample by assumption, and Mori's program for $D$ terminates.
So we can assume it is not nef. Choose a general ample divisor
$A \in \ac  X.$ and look at the intersection point of
the line segment $\overline{AD}$ with the boundary of $\nef X.$.
This defines a $D$-negative contraction. We can assume (by taking
a bigger boundary wall) that it is of relative Picard number one.
If it's small, we can flip it. If it's not birational the program
stops. So we can assume it's a divisorial contraction of relative
Picard number one $f: X \rightarrow Y$. Thus $Y$ is $\bbbq$-factorial.
$D - f^*(f_*(D))$ is effective (since $-D$ is $f$-ample), so we
can replace $D$ by $f^*(f_*(D))$ and assume $D$ is pulled back. Now
we can work in $f^*(\picq Y.) \subset \picq X.$, and
induct on the Picard number of $Y$. Eventually we reduce to the case
when $Y$ has Picard number one, and $D$ is either the pullback
of ample, trivial, or anti-ample. This proves (1).

(4) follows from (1).

Given an effective divisor $D$, running Mori's program for
$D$ yields a birational contraction (indeed a composition of
birational morphisms and flips each of relative Picard number one)
$g: X \dasharrow Y$, with $Y$, $\bbbq$-factorial,
such that $D = g^*(A) + E$ with $A$ semi-ample and $E$ effective
with support the full $g$-exceptional locus. 
Clearly $g^*(A)$ and $E$ are the moving and fixed part of $D$. 
$g^*(\nef Y.) \times \ex g.$ is a Mori chamber by (1.7). This
proves (5).

The contracting morphisms with domain $X_i$ are in one to 
one correspondence with the faces of $\nef X_i.$. Let 
$g: X \dasharrow X'$ be a contracting rational map. Choose
$X_i$ so $g^*(A) \subset f_i^*(\nef X_i)$. It follows that
$X_i \dasharrow X'$ is regular. This proves (2). (3) can
be similarly proved. \qed \enddemo

\remark{1.12 Remark} (1.11.4) is a natural condition, especially
in view of (1.6). 
Unfortunately by itself it does not imply Mori dream space,
or even that nef divisors are semi-ample.
For example let $p: S \rightarrow C$ be the 
projectivization
of the non-split extension of $\ring C.$ by itself, for
$C$ an elliptic curve in characteristic zero. Then the cone of
effective divisors is two dimensional, with edges $F$, the fibre of $p$,
and $C$, the section with trivial normal bundle. Every effective divisor
is nef, and the only non-semi-ample effective  divisor is (a multiple
of) $C$. $\sr S. C.$ is a polynomial ring. Thus all the section rings
are finitely generated. However a natural strengthening of (1.11.4) is
indeed an equivalent characterization of a Mori dream space,
see (2.9). \endremark

\subhead \S 2 Mori Theory and GIT \endsubhead

We refer to \cite{DH} for basic notions from
VGIT. We recall in particular that two $G$-ample line bundles
are called GIT -equivalent if they have the same 
semi-stable locus (and thus in particular give the same
GIT quotients). The equivalence classes are always 
locally polyhedral (and in the cases we consider, will always
be polyhedral). We note one difference from the notation
of \cite{DH}: Here by a GIT chamber we simply mean
a top dimensional GIT equivalence class (in
\cite{DH} the term is reserved for equivalence
classes for which the stable and semi-stable loci are the same). 

\subhead 2.0 Notation \endsubhead
Let $V$ be an affine variety over $k$. Let $G$ be a reductive group
acting on $V$. Let $L$ be the trivial line bundle with the
trivial induced action (i.e. the action is only on the $V$ component).
For each character $\chi \in \chi(G)$, let 
$U_{\chi} = V^{ss}(L_{\chi})$, with quotient 
$q_{\chi}: U_{\chi} \rightarrow U_{\chi}//G \deq Q_{\chi}$. Let
$C \deq C^G(V) \cap ker(f)$ where $f$ is the forgetful map
$f:C^G(V) \rightarrow NS^1(V)$. We denote the complement of the 
semi-stable locus (i.e. the non-semi-stable locus) by $V^{nss}(L_{\chi})$.

\proclaim{2.0.1 Lemma} Notation as above. $C$ is the affine hull of finitely
many characters. \endproclaim
\demo{Proof} This is well known. See e.g., \cite{DH, 1.1.5} or  \cite{T,2.3}. \qed \enddemo

\proclaim{2.1 Lemma} Let $f: U \rightarrow Q$ be a geometric quotient by
a reductive group $G$ acting with finite stabilizers. If $U$ is $\bbbq$-factorial,
and for each $G$-invariant Cartier divisor $D \subset U$, $\ring U.(mD)$ has a linearization
for some $m > 0$, then $Q$ is $\bbbq$-factorial. 

If $G$ is connected then the converse holds.
\endproclaim
\demo{Proof} First we consider the forward implication. Let $D' \subset Q$ be an
effective Weil divisor. Replacing $D'$ by a multiple we may assume the
inverse image, $D$ is Cartier and $\ring U.(D)$ has a linearization. Then 
$D$ is the zero locus of a section, $\sigma$, on which $G$ acts by a character, $\chi$.
Thus if we adjust the linearization, $\sigma$ is an invariant section. The line
bundle and the section descend, by Kempf's descent lemma, after taking
multiples. 

For the reverse direction assume $G$ is connected.
By 
\cite{V,Theorem 1}, since $Q$ is $\bbbq$-factorial, the composition
$$
\picq Q. @>{f^*}>> \picgq U. @>>> \picq U. @>>> A^1(U)_{\bbbq}
$$
is surjective. The first map is an isomorphism by the descent lemma. 
The result follows. \qed \enddemo

\proclaim{2.2 Lemma} Notation as in (2.0). Let $x$ be a character such that
the quotient $Q_x$ is projective. Consider
the following conditions
\roster
\item $V^{ss}(L_x) = V^{s}(L_x)$ and the 
complement $V^{nss}(L_x) \subset V$ has codimension at least two.
\item $V$ has torsion class group.
\item $Q_x$ is $\bbbq$-factorial.
\item Both of the maps 
$$
\chi(G)_{\bbbq} @>{\chi \rightarrow L_{\chi}|_{U_x}}>>  \picgq U_x. 
@<{q_x^*}<< \picq Q_x. 
$$
are isomorphisms.
\endroster
(1-2) imply (3-4). If $G$ is connected, then (1),(3) and (4)
together imply (2). 
\endproclaim
\demo{Proof} Assume (1-2). The the second map in (4) is an isomorphism
by Kempf's descent lemma, and the first map is injective by the codimension
condition of (1). As any two linearizations of a $\bbbq$-line bundle differ 
by a character, (2) implies the first map is surjective.

Assume (1),(3) and (4) and that $G$ is connected. 
$V$ and $U_x$ have the same class group by the
codimension assumption of (1). $U_x$ is $\bbbq$-factorial by (2.1),
and so has torsion class group by the first map in (4). \qed\enddemo

\proclaim{2.3 Theorem} Notation as in (2.0). Let $x$ be a character such
that $Q_x$ is projective. If conditions (2.2.1-2.2.4) hold (e.g. by
the Lemma, if either (2.1.1-2.1.2) hold, or $G$ is connected 
and (2.2.1),(2.2.3), and (2.2.4) hold) then 
$Q_x$ is a Mori dream space. Moreover:
The isomorphism
$\psi: \chi(G)_{\bbbq} \rightarrow \nsd Q_x.$ (induced by (2.1.4) ) identifies 
$\mcd Q_x.$ with $C$, and under this identification, Mori chambers are
identified with GIT chambers. 
Every contraction $f: Q_x \dasharrow Y$ with ($Y$ normal and projective) 
is induced by GIT, i.e. $Y = Q_y$ for some character $y$, and 
$f$ is the induced map.
\endproclaim
\remark{2.3.1 Remark} (2.3) has an obvious analog for quotients of a projective 
variety where we 
vary the linearization on powers of a fixed ample divisor. For
the proof one passes to the cone
over the variety, and applies (2.3).  We leave the details to the reader.
We expect one could further generalize the proposition, to show that GIT
quotients of Mori dream spaces are again Mori dream spaces. 
\endremark

\demo{Proof of 2.3} $\pic Q_x.$ is finitely generated by (2.2.4) and thus
we have (1.10.1). 
Every line bundle on $Q_x$ is of form
$\psi(L_y)$ and $L_y|_{U_x} = q_x^*(\psi(L_y))$. By descent and the codimension
condition we have canonical identifications
$$
H^0(V,L_y)^G = H^0(U_x,L_y)^G = H^0(Q_x,\psi(L_y)). \tag{2.3.3}
$$
Thus $\psi$ identifies $C$ with $\ec Q_x.$.

By the GIT construction,
$L_y|_{U_y}= q_y^*(L_y')$ for an ample line bundle $L_y'$ on $Q_y$, and there
are canonical identifications
$$
H^0(V,L_y)^G = H^0(U_y,L_y)^G = H^0(Q_y,L_y') \tag{2.3.4}
$$
By the codimension condition we have further identifications
$$
H^0(U_y,L_y)^G = H^0(U_y \cap U_x,L_y)^G = H^0(q_x(U_y \cap U_x),\psi(L_y)). \tag{2.3.5}
$$
(Note since $q_x$ is a geometric quotient, $q_x(U_x \cap U_y)$ is open and
its inverse image under $q_x$ is $U_x \cap U_y$).

Every section ring on $Q_x$ is finitely generated by Nagata's
theorem. So Mori equivalence is well defined on the cone of divisors.
Let 
$$
f_{y}: Q_x \dasharrow Q_y
$$
be the induced rational map. 
By (2.3.3) 
$f_{y} = f_{\psi(L_y)}$ and in particular by 
(1.6) a contraction. Further
by (1.6)
$$
\psi(L_y) = f_{y}^*(L_y') + E_y \tag{2.3.6}
$$
for some effective $f_{y}$ exceptional divisor $E_y$. Via $\psi$ we have
both Mori and GIT equivalence on $\mcd Q_x.$.  Clearly GIT equivalence is
finer: 
if the semi-stable loci are the same, the associated contractions
of $Q_x$ are the same. 
By the theory of VGIT, the GIT chambers are finite polyhedral, the affine hulls
of finitely many effective divisors.
Thus $\mcd Q_x.$ is a union of
finitely many Mori chambers, each finite polyhedral.

Now suppose $y$ and $z$ are general members of the same Mori chamber. We will show
they are in the same GIT chamber (thus showing GIT and Mori chambers are the same).
By assumption $f_z$ and $f_y$ are the same, they are birational as the
corresponding divisors are big. By dimension considerations
(since the Mori equivalence class is maximal dimensional), $E_y$ and $E_z$
have the same support, 
the full divisorial exceptional locus of $f_z=f_y$, and the number of components
of either is the relative Picard number and $Q_z = Q_y$ is $\bbbq$-factorial.
We argue now that $U_z = U_y$. 

Of course it is enough to show $U_z \subset U_y$.

Let $z$ be a point of $U_z$. Then by the construction of GIT quotients,
there is a section $\sigma \in H^0(V,L_y)^G$ such that
$\sigma|_{U_y} = q_y^*(\sigma')$ for a section 
$$
\sigma' \in H^0(Q_y,L_y')
$$
which does not vanish at $q_z(z) \in Q_z = Q_y$. We claim that 
$$
L_y|_{U_z} = q_z^*(L_y') \text{ and } \sigma|_{U_z} = q_z^*(\sigma'). \tag{2.3.7}
$$
This will of course imply
$\sigma(z) \neq 0$ and $z \in U_y$. We can check (2.3.7) after removing
any codimension two subset from $U_z$. 
By (2.3.3) and (2.3.6) $U_x \cap U_y$ and $U_x \cap U_z$ are equal in codimension one: the
complement of either is, up to codimension one, 
the inverse image under $q_x$ of the divisorial exceptional
locus of $f_z = f_y$. Thus $U_z$ and $U_y$ are equal in codimension one, and we
can check (2.3.7) after restricting to $U_z \cap U_y$ where it obviously holds.

Thus the Mori and GIT chambers have the same interiors, and (up to
closure) each chamber is of form $f_z^*(\ac Q_z.) \times \ex f_z.$ 
for linearizations $z$ such that $Q_z$ is $\bbbq$-factorial.
In particular (up to closure) the moving cone will be the 
union of the (finitely many) chambers with $f_z$ small. To finish the proof
we need only show that on these $Q_z$ the nef cones are generated by
finitely many semi-ample line bundles. Let $z$ be such a character.
Note since $f_z$ is small, $\mcd Q_z.$ and $\mcd Q_x$ are canonically
identified by $f_z^*$. Let
$C_z \subset C$ be the closure of the GIT chamber of $z$ (which we know is
the closure of the ample cone of $Q_z$).
Choose $y \in \partial C_z$. By the VGIT theory
there is an inclusion
$U_z \subset U_y$. It follows that 
the rational map 
$$
f_{zy}= f_y \circ f_z^{-1}: Q_z \dasharrow Q_y
$$
is regular. By negativity
of contraction, since $\psi(L_y)$ (being on the boundary of the ample 
cone) is nef on $Q_z$ the term $E_y$ in (2.3.6) is empty, and 
$\phi(L_y) = f_{y}^*(L'_y)$. Since $L'_y$ is ample, and $f_{zy}$ is regular,
$\psi(L_y)$ is semi-ample
on $Q_z$. 
\qed \enddemo

\proclaim{2.4 Corollary} Let $X$ be a projective geometric GIT quotient for 
the action of an algebraic torus on an affine variety with torsion
class group. If the non-stable locus has codimension at least two, then 
$X$ is a Mori dream space satisfying the conclusions of (2.3).
Moreover
GIT quotients from linearizations in the interior of
Mori chambers are geometric quotients (i.e. the Mori chambers are
chambers in the sense of \cite{DH}.)

Suppose furthermore that $V$ is smooth. Then any rational contraction
of $X$ with $\bbbq$-factorial image 
is a composition of weighted flips, weighted blowdowns,
and \'etale locally trivial (on the image) fibrations of relative Picard
number one with fibre a quotient of weighted projective space 
by a finite abelian group (in particular the image
of any such a contraction has cyclic quotient singularities). Indeed the
factorization is obtained by the series of (necessarily codimension
one) wall crossings connecting a general member of the ample cone of $X$
with a general member of the chamber corresponding to the contraction.
\endproclaim

The smooth case of (2.4) is obviously an optimal situation --the 
contractions are parameterized in a nice combinatorial way, and 
each contraction is naturally factored into simple parts. We note
that in general such a factorization is only possible once one
allows small modifications --there will be no such factorization if one restricts
themselves to morphisms. e.g. there are birational morphisms
$f:X \rightarrow Y$ of relative Picard number two 
between smooth projective toric varieties which do not factor through
any morphism $X \rightarrow Y'$ of relative Picard number one with $Y'$
$\bbbq$-factorial. 

\demo{Proof of (2.4)} Except for the final claim of the first paragraph, everything is
immediate from (2.3) and the theory of VGIT 
\cite{DH,0.2.5} or \cite{T,5.6}. We follow the notation of the proof of (2.3).
Consider a linearization $y$ in
the interior of a Mori chamber. 
It's enough to show $q_y^*$ is an isomorphism:
Then for any character $v$, $L_{mv}|_{U_y}$ is pulled back from $Q_y$ (for some $m>0$), thus
the stabilizer of any point of $U_y$ is in the kernel of $mv$ for all $v$, and so 
the stabilizer is finite. We can check $q_y^*$ is an isomorphism after removing
codimension two subsets from $Q_y$ and $U_y$. Thus we can restrict to
$U_y \cap U_x$ and to the locus where $f_y^{-1}$ is an isomorphism. Here 
the quotient is geometric, so $q_y^*$ is an isomorphism by Kempf's descent
lemma. \qed \enddemo
 

(2.4) applies to any quasi-smooth projective toric variety $X$ by
Cox's construction which gives an essentially canonical
way of writing $X = X(\Delta)$ (for the fan $\Delta$ with support the lattice
$N = \bbbn^n$) 
as a GIT quotient of $\bbba^r$,
$r = \#(\Delta(1))$ (where $\Delta(k)$  is the collection of $k$-dimensional
cones in the fan).
by $T = \Hom A_{n-1}. \gm .$, satisfying
(2.1.1-2). See \cite{C}.

For a $\rho$ dimensional torus $T$ acting on affine space the GIT chambers
are particularly simple: An action of $T$ on $\bbba^r$  is given by 
$r$ characters $\chi_i \in \chi(T)$. The $T$-ample cone is the affine hull
of the characters. The GIT chambers are the affine hulls of all subsets
of $\rho$ independent characters. See e.g., \cite{DH}. Combining this
with Cox's construction and (3.3) gives a simple algorithm for describing
the Mori chambers of any quasi-smooth projective toric variety. This description
was obtained by Oda and Park,  \cite{OP}, using Reid's Toric Mori's program 
\cite{R1}. The factorization in (2.4) gives a cheap form of Morelli's
factorization theorem, \cite{Mo}, cheap in that even in factoring
a birational map between smooth spaces we allow cyclic quotient
singularities. The factorization does however have an important advantage over
Morelli's: Morelli factors birational maps, but even to factor a birational
morphism he may have to blowup an indeterminate number of times. In
particular there could \'a prior be infinitely many such factorizations.
On the other hand all possible factorizations into contractions are encoded
in the chamber decomposition of (2.4). 

We note that by \cite{BK} quasi-smooth projective spherical varieties
give further examples of Mori dream spaces.

Notation: For a collection of $r$ line bundles $L_1,\dots,L_r$ and 
a vector of integers $v=(n_1,\dots,n_r) \in \bbbz^r$ we use the notion
$$
L^v \deq L_1^{\otimes n_1} \otimes L_2^{\otimes n_2} \dots \otimes L_r^{\otimes n_r}.
$$

\definition{2.5 Definition} For line bundles $L_1,\dots,L_r$ on $X$ let
$$
\sr X. L_1,\dots,L_r. \deq \bigoplus_{v \in \bbbn^r} H^0(X,L^{v}).
$$

\enddefinition

\definition{2.6 Definition} Let $X$ be a projective variety such
that $\picq X. = \nd X.$.  By a {\it Cox ring} for $X$ we mean the
ring 
$$
\cox X. \deq \sr X. L_1,\dots,L_r.
$$
for a choice of line bundles $L_1,\dots,L_r$ which are a
basis of $\picq X.$, and whose affine hull contains $\ec X.$. 
\enddefinition

\remark{Remark}
Rather than have choices as in (2.6) we would prefer to use 
$$
\bigoplus_{L \in \pic X.} H^0(X,L)
$$
but this does not have a well defined ring structure (for an isomophism
class $L$ the vector space
$H^0(X,L)$ is only determined up to a scalar). 
Of course $\cox X.$ as we've defined it depends on the choice of basis. If 
we choose two $\bbbz$-basis of the torsion free part of $\pic X.$ then
the two rings are isomorphic. If we replace the line bundles by positive powers,
the original Cox ring is a finite extension of the new Cox ring. Thus 
finite generation of the Cox ring, which for our purposes will be the
main issue, is independent of choice. For any toric variety
$\cox X.$ is a polynomial ring, Cox's coordinate ring, \cite{C}, whence
the name. \endremark

\proclaim{2.7 Lemma} Let $\sigma_1,\sigma_2 \in H^0(X,L)$ be two
sections of a non-torsion line bundle whose zero divisors have no
common component. Then $(\sigma_1,\sigma_2) \subset \cox X.$ is 
a regular sequence.
\endproclaim
\demo{Proof} Suppose $a \cdot \sigma_1 = b \cdot \sigma_2$. 
We can assume $a$ and $b$ are homogeneous. Thus
$a$ and $b$ are sections of the same line bundle $M$ and 
$a \otimes \sigma_1 = b \otimes \sigma_2$. Let $A,B,D_1,D_2$ be the
zero divisors of $a,b,\sigma_1,\sigma_2$. We have 
an equality of Weil divisors
$$
A + D_1 = B + D_2
$$
It follows that $A - D_2 = B -D_1$ is effective and Cartier. Thus
$d = a/\sigma_2 = b/\sigma_1$ is a regular section of $M \otimes L^{*}$,
and $a = \sigma_2 \cdot d$, $b = \sigma_1 \cdot d$. \qed \enddemo

\proclaim{2.8 Lemma (Zariski)} Let $L_1,\dots,L_d$ be semi-ample line bundles
on a projective variety $Y$. Then $\sr Y. L_1,\dots,L_d.$ is 
finitely generated, and there exists an integer $ m >0$, such that
for any $k > 0$, 
after replacing $L_i$ by $L_i^{\otimes km}$ the canonical map
$$
H^0(Y,L_1)^{\otimes n_1}\otimes \dots H^0(Y,L_d)^{\otimes n_d} \rightarrow 
H^0(Y,L^{(n_1,\dots,n_d)})
$$
is surjective for all $n_i \geq 0$.
\endproclaim
\demo{Proof} If $\bbbp = \bbbp(L_1 \oplus \dots \oplus L_r)$ then 
$\sr Y. L_1,\dots,L_r. = \sr {\bbbp}. {\ring {\bbbp}.(1)}.$, so we 
we reduce to a single 
semi-ample line bundle, where finite generation is a familiar result
due to Zariski. For some $m > 0$, 
$$
\sym k. {H^0(\bbbp,\ring.(m)) }. \rightarrow H^0(\bbbp,\ring.(km) )
$$
is surjective for all $k$. The second statement follows by considering
the appropriate graded piece. \qed \enddemo

\proclaim{2.9 Proposition} Let $X$ be a $\bbbq$-factorial projective
variety such that $\picq X. = \nd X.$. 
$X$ is a Mori dream space iff 
$\cox X.$ is finitely generated. 

If $X$ is a Mori dream space then $X$ is a GIT quotient of
$V = \spec {\cox X.}.$ by the torus 
$G = \Hom {\bbbn}^r. \gm .$, where $r$ is the Picard number of $X$,
satisfying the conditions of (2.3).
Moreover we may choose the Cox ring so that 
$G$ acts freely on the semi-stable loci of any linearization in the interior
of a Mori chamber.
\endproclaim
\demo{Proof} 
Let $R = \cox X. = \bigoplus_{v \in \bbbn^r} R_v$.

Assume $X$ is a Mori dream space. 
For each (closed) Mori chamber $C \subset \ec X.$
let $R_C = \bigoplus_{v \in C} R_v$. As there are only finitely many chambers, and any 
homogenous element of $R$ lies in some $R_{C}$, to show $R$ is finitely generated it 
is enough to show $R_C$ is finitely generated for each $C$. Choose a chamber $C$ 
and line bundles $J_1,\dots, J_d \in C$ which generate $C$ (as a semi-group). Expressing
the $J_i$ as tensor products of the $L_j$ induces a 
surjection 
$$
\sr X. J_1,\dots,J_d. \rightarrow R_C
$$
so we need only show $\sr X. J_1,\dots,J_d.$ is finitely generated.
By
(1.11.2) there is a 
contracting rational map $f:X \dasharrow Y$ to a projective $\bbbq$-factorial
normal variety $Y$ such that each $J_i = f^*(A_i)(E_i)$ for a $A_i$ semi-ample,
and $E_i$ effective and $f$ exceptional. Thus by the projection formula there is 
a natural identification
$$
\sr X. J_1,\dots,J_d. = \sr Y. A_1,\dots,A_d.
$$
The latter is finitely generated by (2.8).

Now suppose $R$ is finitely generated. Note $G$ acts naturally on $R$ so that
$$
R = \bigoplus_{v \in \chi(T)=\bbbn^r} R_v
$$
is the eigenspace decomposition for the action. Thus 
$$
H^0(V,L_v)^G = R_v
$$
and for $v = L \in \pic X.$ the ring of invariants is
$$
{\sr V. L_v.}^G = \sr X. L. .
$$
Thus $X$ is the GIT quotient for any linearization $v \in \ac X. \subset \chi(G)_{\bbbq}$,
and for any linearization $v = L$ the induced rational map $X \dasharrow Q_v$ is
$f_L$.

Let $h: R \rightarrow \bbbc$ be a point of $V$, and $v \in \ec X.$ a
linearization. By our description of
the invariants $h$ is $L_v$ semi-stable
iff $h(R_{nv}) \neq 0$ for some $n > 0$. Suppose for some $m>0$ and
$v_1,\dots,v_d$ with $\sum v_i = mv$ that
$$
R_{nv_1} \otimes \dots \otimes R_{nv_d}
\rightarrow R_{nmv}
$$
is surjective for all $n > 0$. Then 
$$
V^{ss}(L_v) = \bigcap_{i =1}^{i=d} V^{ss}(L_{v_i}).
$$ 
It follows in particular from (2.8) that
any ample $v$ has the same semi-stable locus, say $U$. Furthermore, $\lambda \in G_h$ (the 
stablizer of $h$) iff $\lambda^{v} = 1$ for all $v$ such that $h(R_v) \neq 0$. In
particular $\lambda^{v}$ is torsion if $h$ is $L_v$ semi-stable. The ample cone
generates $\bbbn^r = \chi(T)$ (as a group) thus any $h$ semi-stable for an ample $v$ has
finite stabilizer. Thus $X$ is a geometric quotient of $V$. Choose two
sections $\sigma_1,\sigma_2$ of some ample line bundle $L$, whose zero divisors
have no common component. Let $I$ be the ideal of the non semi-stable
locus $U^{c}$ (with reduced structure).
$\sigma_1,\sigma_2 \in I$, so by (2.7), $U^c$ has codimension at least two. 

Thus the quotient $X$ satisfies the conditions of (2.3), so $X$ is a Mori
dream space.

Now choose a Mori chamber $C$ and generating line bundles $J_1,\dots,J_d$,
with associated contracting birational map $f$ as above. 
After replacing $X$ by a SQM (which is again a Mori dream
space with the same Cox ring) we may assume $f$ is a morphism. One sees 
that $V^{ss}(L_v)$ is constant for $v$ in the interior of $C$, and that any point
in this open set has finite stablizers using (2.8) exactly
as in the case of $C = \nef X.$ argued above.
The same argument shows that after replacing the $L_i$ by powers, the stabilizers
are trivial. \qed \enddemo

\proclaim{2.10 Corollary} Let $X$ be a smooth  projective
variety with $\picq X. = \nd X.$. $X$ is a toric variety iff it has
a Cox ring which is a polynomial ring. 
\endproclaim
\demo{Proof} In the smooth toric case $\cox X.$ is Cox's homogeneous coordinate
ring. By (2.9) if  $\cox X.$ is finitely generated, then $X$
is a geometric GIT quotient of $\spec {\cox X.}.$ by a torus, and 
and the quotient of affine space
by a torus is a toric variety. \qed \enddemo

The next proposition indicates that the birational contractions of 
a Mori dream space are induced from toric geometry.

\proclaim{2.11 Proposition} Let $X$ be a Mori dream space. Then there is embedding
$X \subset W$ into a quasi-smooth projective toric variety such that 
\roster
\item The restriction $\picq W. \rightarrow \picq X.$ is an isomorphism
\item The isomorphism of (1) induces an isomorphism $\ec W. \rightarrow \ec X.$.
\item Every Mori chamber of $X$ is a union of finitely many Mori chambers of $W$.
\item For every rational contraction $f: X \dasharrow X'$ there is toric
rational contraction $\tilde{f}: W \dasharrow X'$, regular at the generic
point of $X$, such that $f = \tilde{f}|_X$.
\endroster
\endproclaim
\demo{Proof} 
Let $R = \cox X. = \bigoplus_{v \in \bbbn^r} R_v$. By (2.9) $R$ is finitely
generated over $R_0 = k$. Choose homogenous generators whose degrees 
(in the grading) are non-trivial effective divisors. This defines
a $k$-algebra surjection $A \rightarrow R$ from a polynomial ring $A$,
and a compatible action of $T = \Hom {\bbbn}^r. {\gm}.$ on $A$
such that $A^T = k$. Let $\bbba = \spec A.$. We have an equivariant
embedding
$V = \spec R. \subset \bbba$. Let $M_v$ (resp. $L_v$) be twistings by the character
$v \in \chi(T)$ of the trivial line bundle  on $\bbba$ (resp. $V$). 
We follow
the notation of (2.0).
$\bbba^{ss}(M_v) \cap V = V^{ss}(L_v)$ for any $v$.  Thus GIT equivalence
on $V$ is finer than GIT equivalence on $\bbba$. 
Choose $v$ a general member of a Mori chamber of $\mov X.$. We claim that
the quotient $W_v \deq \bbba^{ss}(M_v)//T$ satisfies the conditions
of (2.3). As
remarked in the proof of (2.4) we need only check the codimension condition
of (2.2.1). Suppose
${\bbba}^{nss}(M_v)$ has a divisorial component. By (2.9) $Q_v$ satisfies
the conditions of (2.3) , so 
there is a non-constant function
$g \in \ring.(\bbba)$ on which $T$ acts by some character, $\chi$, whose
restriction to $V$ is a unit. But then $L_{\chi} \in \picg U_v.$ is trivial.
So $\chi$ is trivial. But then $f$ is a non-constant invariant function,
a contradiction. Thus (2.3) applies to $Q_v$ and $W_v$. The result follows.
\qed \enddemo
 
There is a natural local (in the cone of divisors) generalization
of (1.10).

\definition{2.12  Definition} Let $C \subset \ec X.$ be the affine hull
of finitely many effective divisors. We say that $C$ is a {\bf Mori 
Dream Region} provided the following hold:
\roster
\item There exist a finite collection of birational contractions
$f_i:X \dasharrow Y_i$ such that $C_i \deq C \cap f^*(\nef Y_i.) \times \ex{f_i}.$
is the affine hull of finitely many effective divisors.
\item $C$ is the union of the $C_i$.
\item Any line bundle in $(f_i)_*(C_i) \cap \nef Y_i.$ is semi-ample.
\endroster
\enddefinition

(2.9) has the following analog
\proclaim{2.13 Theorem} Let $X$ be a normal projective
variety and let $C \subset \nd X.$ be a rational polyhedral
cone (ie the affine hull of the classes of finitely many line bundles).
$C \cap \ec X.$ is a Mori dream region iff there are generators
$L_1,\dots,L_r$ of $C$ such that $\sr X. L_1,\dots,L_r.$ is 
finitely generated.
\endproclaim
\demo{Proof} Analogous to that of (2.12). \qed \enddemo

It is natural to expect that the 
region of the cone of divisors studied by Mori theory is itself (at least locally)
a Mori dream region. This leads to the following conjecture, which by the
ideas of the proof of (2.9) contains all the main conjectures/theorems 
(e.g. cone and contraction
theorems, existence of log flips, and log abundance) of Mori's program:

\proclaim{2.14 Conjecture} Let $\Delta_1,\dots,\Delta_r$ be a collection
of boundaries such that $K_X + \Delta_i$ is Kawamata log terminal. Choose an integer
$m$ so that $L_i =m(K_X + \Delta_i)$ are all Cartier. Then 
$\sr X. L_1,\dots,L_r.$ is finitely generated. 
\endproclaim

\proclaim{2.15 Corollary } The conjecture holds in dimension three or
less. \endproclaim
\demo{Proof} By \cite{Sh,6.20} the intersection of the affine
hull of the $L_i$ with $\ec X.$ is a Mori dream region.
\qed \enddemo

\proclaim{2.16 Corollary} Let $X$ be a log Fano $n$-fold, with $n \leq 3$. Then $X$
is a Mori dream space. \endproclaim
\demo{Proof} Let $K_X + \Delta$ be KLT and anti-ample. 
Choose a basis $L_1,\dots,L_r$ of $\pic X.$ whose affine
hull contains $\ec X.$. Choose $n>0$ so that $A_i = L_i - n(K_X + \Delta)$
is ample for all $i$. Let $\Delta_i = 1/nm D_i + \Delta$ for $D_i$ a general
member of $|mA_i|$. Note $L_i = n(K_X + \Delta_i)$ (in $\picq X.$), and
$\Delta_i$ is KLT for sufficiently large $m$.  Now apply (2.13).
\qed \enddemo

\subhead \S 3 Connections with $\mon$ \endsubhead

The original motivation for this paper was to try to understand
the geometric meaning of the cone of effective divisors, 
in connection with questions of
Fulton on $\mon$,
the moduli space of stable $n$-pointed rational curves. 

\proclaim{3.1 Question (Fulton)} Is $\mc {\mon}.$ (resp $\ec {\mon}.$) 
the affine hull of the one dimensional (resp. codimension one) strata?
\endproclaim

See \cite{KeM} for definitions and partial results, and an indication of
the wide range of contexts in which $\mon$ naturally appears. 

The connection with GIT is as follows:

Consider the diagonal action of $G=\Sl 2.$ on the $n$-fold product 
$(\pr 1.)^{\times n}$. By the
Gelfand-Macpherson correspondence, the VGIT theory for this action is identified
with that of the torus $T = {\Bbb G}_m^n$ on the Grassmannian $G(2,m)$, e.g. the $G$-ample
cones and their chamber decompositions are naturally identified, and the corresponding
GIT quotients are the same (in the first case we vary the line bundle, the
linearization on each is canonical, in the second the line bundle is fixed and we
vary the linearization by characters). (2.4) (see (2.3.1)) now applies.
The $G$-ample cone and chamber decomposition
are easy to describe, see \cite{DH}, and one obtains a complete description of
the rational contractions on any of the GIT quotients. 
By \cite{K}, $\mon$ is the inverse limit of all the
GIT quotients. 

\proclaim{3.2 Question} Is $\mon$ a Mori dream space?
\endproclaim

One result of \cite{KeM}
is that any extremal ray
of $\mc {\mon}.$ which can be contracted by a map of relative Picard number
one is generated by a stratum, 
so long as the exceptional locus of the map has dimension at least two
(any stratum can be contracted, and the exceptional locus of the 
contraction satifies the dimension condition for any $n \geq 9$).
By (1.11)
if $\mon$ is a Mori dream space, then any extremal ray of the Mori cone 
is contracted by a map
of relative Picard number one. Thus a positive answer to (3.2) would
nearly answer Fulton's question
for $\mc {\mon.}$. 

There is a natural action of the symmetric group $S_n$ on $\mon$ and it
is natural to consider the $S_n$ equivariant geometry, or equivalently
the geometry of the quotient $\tmon$. This quotient is itself an important moduli
space, e.g. $\widetilde{M}_{0,2g + 2} \subset \mg$ is the hyperelliptic
locus.

Let $k = [n/2]$. Let $B_i \subset \mon$,
$k \geq i \geq 2$, be the union of codimension one strata whose
generic point corresponds to a curve with two components, and exactly $i$ 
marked points on one of the components. 
The analog of (3.1) for $\ec {\tmon}.$ is proven in \cite{KeM}.
$\ec {\tmon}.$ is in fact simplicial, generated by the (images of the) $B_i$.
Furthermore, every moving divisor is big (thus every rational contraction
of $\tmon$ is birational).

In particular by (2.9), $\tmon$ is  a Mori dream space iff ring 
$$
\bigoplus_{(d_2,\dots,d_k) \in {\bbbn}^{k-2}} H^0(\mon,\sum d_i B_i)
$$
is finitely generated.

Furthermore, by (1.10), a positive answer to (3.2) would imply the
following:

\proclaim{3.4 Implication} For each $ k \geq i \geq 2$ there exists 
a unique birational contraction $f_i: \tmon \dasharrow Q_i$, 
where $Q_i$ is $\bbbq$-factorial of Picard number one,
and the exceptional divisors of $f_i$ are exactly the $B_i$
with $j \neq i$. The moving cone of $\tmon$ is simplicial, generated
by pullbacks of ample classes from the $Q_i$. 
\endproclaim

$f_2$ of (3.4) exists, it is the (regular) contraction to the 
GIT quotient of $\Sl 2.$ for the action on the $n^{th}$ symmetric
product of the standard representation (i.e. the GIT quotient
for $n$ unmarked points on $\pr 1.$).

We finish by giving a result which gives another connection between
$\mon$ and GIT. Though rather unrelated to the rest of the paper, we
hope the reader will find it of interest.

In \cite{FM} Fulton and MacPherson construct a functorial 
compactification $\fm X. n.$ of the locus of distinct points in
a smooth variety $X$. As we now indicate, $\mon$ occurs as a GIT
quotient of $\fm {\pr 1.}. n.$ by the natural action of 
$G = \Sl 2.$. 

There is a proper birational morphism $f: \fpone \rightarrow \ponen$.
Let $E$ be an effective divisor, with support the full exceptional
locus of $f$, such that $-E$ is $f$-ample (such an $E$ exists for
any proper birational morphism between $\bbbq$-factorial varieties).

\proclaim{3.5 Theorem} For each linearization $L \in \picG \ponen.$ such
that $(\ponen)^{ss}(L) = (\ponen)^{s}(L) \neq \emptyset$, and each
sufficiently small $\epsilon > 0$ the line bundle
$L' = f^*(L)(-\epsilon E)$ is ample and 
$$
(\fpone)^{ss}(L') = (\fpone)^{s}(L') = f^{-1}(\ponen)^{ss}(L).
$$
There is a canonical identification
$$
(\fpone)^{ss}(L') / G = \mon 
$$
and a commutative diagram
$$
\CD
(\fpone)^{ss}(L') @>{f}>> (\ponen)^{ss}(L) \\
@V{q}VV                       @V{q}VV \\
\mon  @>g>>                   (\ponen)^{ss}(L)/G
\endCD
$$
where $q$ indicates the geometric quotient, and the various $g$
are Kapranov's blowup expressions for $\mon$, realizing it as the
inverse limit of all the GIT quotients of $\ponen$.
\endproclaim
\demo{Proof} We follow the notation of \cite{FM} for divisors on
$\fm X. n.$ and \cite{H1} for the chamber decomposition for 
$\picG {\ponen}.$. For a subset $S \subset \{ 1,2,\dots,n\}$
let $l_S$ be the linear functional on $\picG {\ponen}.$
$$
l_S(x_1,\dots,x_n) =  \sum_{i \in S} x_i - \sum_{i \not \in S} x_i.
$$

For the first statement see \cite{H2}.
We let $U$ be the semi-stable locus for a linearization
on $\ponen$ corresponding to a chamber, and $U' = f^{-1}(U)$. Let
the corresponding quotients be $Q$ and $Q'$. By \cite{FM,pg. 195, pg. 212}, there
is a natural $G$-equivariant surjection $\fpone \rightarrow \mon$,
where $G$ acts trivially on $\mon$. 
Thus there is an induced proper birational morphism $Q' \rightarrow \mon$.
To prove this is an isomorphism, both sides being $\bbbq$-factorial,
it is enough to show that both sides have the same Picard number.
$$
\rho(Q') = \rho(U') = \rho(U) + e_U = \rho(Q) + e_U
$$ 
where $e_U$ is the number
of $f$-exceptional divisors that meet $U'$, or equivalently,
the number of diagonals $\Delta_S$ which meet $U$ and have $|S| > 2$.
We show first that $\rho(Q')$ is constant (i.e. independent of the 
chamber). It's enough to check two
chambers sharing the codimension one wall $W_S$. Let the
two open sets be $U_1, U_2$ where we assume $U_1$ meets $\Delta_S$
and $|S| \leq |S^c|$. Note the $U_i'$ meet the same divisors
$D(T)$ except that $U_1'$ meets $D(S)$ (and not $D(S^c)$)
while $U_2'$ meets $D(S^c)$ (and not $D(S)$). 
If $|S| \geq 2$ then $Q_1 \dasharrow Q_2$ is
a small modification, so $\rho(Q_1) = \rho(Q_2)$, and $e_{U_1} =
e_{U_2}$.
Suppose $|S| = 2$. Then $Q_1 \dasharrow Q_2$ is a birational
contraction with
exceptional divisor (the image of) $\Delta_{S}$. Thus
$\rho(Q_1) = \rho(Q_2) + 1$. On the other hand $e_{U_1} = e_{U_2} -1$,
since $D(S)$ is not exceptional (its image is divisorial), while
$D(S^c)$ is exceptional. 

Now we compute $\rho(Q')$ for the case of the chamber given by
ineqalities $l_{S} < 0$ for all $1 \not \in S$.
In this case $Q = \pr {n-3}.$, while the $f$ exceptional
divisors that meet $U'$ are precisely the $D(S)$ with 
$1 \not \in S$, and $ n-2 \geq |S| \geq 3$. Thus 
$$
\rho(Q') = 2^{n-1} - {\binom n 2} -1 = \rho(\mon). \qed
$$ 
\enddemo

\medskip
\ref \by [BK] \quad M. ~Brion and F. ~Knop
      \paper Contractions and flips for varieties with group action of small complexity
      \jour  J. Math. Sci. Univ. Tokyo 
      \vol 1 
      \yr 1994
      \pages 641--655
\endref
\ref \by [BP] \quad M. ~Brion and C. ~Procesi
     \paper Action d'un tore dans une vari\'et\'e projective
      \jour Progress in Math.
      \vol 192
       \year 1990
       \pages 509--539
\endref
\ref \by [C] \quad D. ~Cox
      \paper The homogeneous coordinate ring of a toric variety
      \vol 4
      \jour Jour. Algebraic Geom.
      \year 1995
      \pages 17--50
\endref
\ref \by [DN] \quad J. ~Drezet and M. ~Narasimhan
     \paper Groupe de Picard des vari\'et\'es de modules de fibr\'es
semistable sur les courbes alg\'ebriques
     \jour Inven. Math.
     \vol 97
      \yr 1989
      \pages 53--94
\endref
\ref \by [DH] \quad I. ~Dolgachev and Y.~Hu
     \paper Variation of geometric invariant theory quotients
     \jour Inst. Hautes \'Etudes Sci. Publ. Math.
     \vol 87
     \yr 1998
     \pages 5--56
\endref
\ref \by [FM] \quad W. ~ Fulton and R. ~MacPherson,
    \paper A compactification of configurations spaces
     \jour Ann. of Math.
     \vol 139
      \year 1994
      \pages 183--225
\endref
\ref\by [H] \quad Y. ~Hu
    \paper The Geometry and Topology of Quotient Varieties of Torus Actions
     \jour Duke Math. Jour.
     \vol  68
     \year 1992
     \pages 151-183   
\endref
\ref\by [H1] \quad Y. ~Hu
     \paper Moduli of stable polygons and symplectic structures on $\mon$
      \jour Compositio Mathematicae 
      \vol 118  
      \year 1999
       \pages 159-187        
\endref
\ref \by [H2] \quad Y. ~Hu 
     \paper Relative Geometric Invariant Theory and Universal Moduli Spaces
     \jour International Journal of Mathematics  
     \vol  7 
     \year 1996
     \pages  151 -- 181.
\endref
\ref \by [K] \quad M. M. Kapranov
     \paper Chow Quotients of Grassmannians. I.
     \inbook I.M. Gelfand Seminar, S. Gelfand, S. Gindikin eds. 
     \yr 1993
     \publ  A.M.S. \bookinfo Advances in Soviet Mathematics vol. 16, part 2.
     \pages 29-110
\endref      
\ref \by [Ka] \quad Y. ~Kawamata
     \paper On the cone of divisors of Calabi-Yau fibre spaces
      \jour Internal J. Math.
      \vol 8
      \yr 1997
      \pages 665-687
\endref
\ref \by [Ke] \quad S. ~Keel
     \paper Basepoint freeness for nef and big linebundles in positive characteristic
     \jour Annals of Math.
     \yr 1999
      \pages 253--286
\endref
\ref\by [KeM] \quad S. ~Keel and \quad J. ~McKernan
    \paper Contractible extremal rays of $\overline{M}_{0,n}$
    \jour preprint alg-geom/9707016
\endref
\ref\by [Koll\'aretal] \quad J.~Koll\'ar (with 14 coauthors)
    \paper Flips and Abundance for Algebraic Threefolds
    \jour Ast\'erique
    \yr 1992
    \vol 211
\endref
\ref\by [KM] \quad J. ~Koll\'ar and S. ~Mori
     \book Birational Geometry of Algebraic Varieties
    \publ Cambridge Univ. Press
    \vol 134
    \yr 1998
\endref
\ref \by [Mo] \quad R. ~Morelli
     \paper The birational geometry of toric varieties
     \jour  J. Algebraic Geom.
     \vol  5 
      \yr 1996 
      \pages 751--782
\endref
\ref \by [M] \quad D. ~Morrison
     \paper Comactifications of moduli spaces inspired by mirror symmetry
      \jour Asterique
      \vol 218
      \year 1993
      \pages 243--271
\endref
\ref \by [OP] \quad T. ~Oda and H. ~Park
      \paper Linear Gale transforms and Gelfand-Kapranov-Zelevinskij  decompositions
      \jour Tohoku Math. Jour. 
       \vol 43
       \year 1991
       \pages 375--379
\endref
\ref \by [R1] M. ~Reid
    \paper Decomposition of toric morphisms. 
     \jour Progr. Math.
     \vol. 36
     \year 1983
     \pages  395--418
\endref
\ref \by [R2] \quad M. ~Reid
     \paper What is a flip?
     \yr 1992
     \jour preprint
\endref
\ref \by [T] \quad M.~Thaddeus
    \paper Geometric Invariant Theory and Flips
    \yr 1996
    \jour Journal of the A. M. S.
    \vol 9
    \pages 691--723
\endref
\ref \by[V] \quad A. ~Vistoli
     \paper Chow groups of quotient varieties
      \jour Jour. of Algebra
      \year 1987
      \vol 107
      \pages 410--424
\endref
\end